\documentclass[11pt]{article}
\usepackage{amsmath, amssymb}
\usepackage{cite}

\newtheorem{theorem}{Theorem}
\newtheorem{theorem*}{Theorem}

\newtheorem{lemma}{Lemma}
\newtheorem{proposition}{Proposition}
\newtheorem{definition}{Definition}
\newtheorem{definition*}{Definition}
\newtheorem{remark}{Remark}

\newtheorem{problem}{Problem}

\begin{document}

\title{\LARGE \bf
 Braids, links, cobordisms and formal groups}
\author{ N.M. Glazunov}

\date{}

\maketitle
\begin{center}
{\small
{\rm  Email:} {\it glanm@yahoo.com }}
\end{center} 



\abstract{
V.V. Sharko in his papers and books has investigated functions on manifolds and cobordism. Braids intimately connect with functions on manifolds. These connections are represented by mapping class groups of corresponding discs, by fundamental groups of corresponding punctured discs, and by some other topological or algebraic structures.
This paper presents selected algebraic  methods and results of braids, links, cobordism connect with investigations by V.V. Sharko. These includes group theoretic results on braids and links, infinitesimal braid group relations and  connections as well as connections on coherent sheaves on smooth  schemes, a sketch of  our algorithm for constructing of Lazard`s one dimensional universal commutative formal group and selected  results on applications of commutative formal groups to cobordism theory.}



\section*{Introduction and Acknowledgments}

V.V. Sharko in his book \cite{Sharko1993} and in papers  has investigated functions on manifolds and cobordisms \cite{Sharko1991,SharkoK1991}. 
Braids intimately connect with functions on manifolds.
These connections are represented by mapping class groups of corresponding discs, by fundamental groups of corresponding punctured discs, and by some other topological or algebraic structures.
Cobordisms and h-cobordisms were studied by V.V. Sharko in connection with the development of Morse theory.
This paper is organized as follows. Section 1 and 2 contain preliminaries
on braids, lincs and their relations with functions on manifolds.
Section 3 in the framework present very shortly results by  Reidemeister, Alexander and Markov.
Sections 4 and 5 discuss infinitesimal braid group relations and  connections as well as connections on coherent sheaves on smooth  schemes. Sections 6 and 7 describe selected results on cobordism theory. 
The author's acquaintance with braids and links, as well as with the theory of cobordism, took 
place in 1970-1971. during his internship with academician A.A. Markov. (A.A. Markov was the head 
of the laboratory at the Computing Center of the USSR Academy of Sciences, and, at the same
 time, the head of the department of mathematical logic at Moscow State University).
  A.A. Markov gave the description of the set of isotopy classes of oriented links in ${\mathbb R}^3$ in terms of braids. For manifolds of the dimension grater than 3 A.A. Markov  has proved the undecidability  of the problem of homeomorphy.
The author first got acquainted with the elements of the theory of cobordism at the organized by S.P. Novikov School "Topology and Algebraic K-Theory", which took place in Pushchino on the Oka in 1971. Academician A. Markov recommended the author to participate in the School.
V.V. Sharko  supported the author's research on the use of formal groups in topology  and algebra.
 I am deeply grateful to them, as well as to the organizations in which I was, for the hospitality and support.
More recently author got acquainted with approaches to the algebraic theory of cobordism at the school "Axiomatic, enriched and motivic homotopy theory'', Isaac Newton Institute for Mathematical Sciences, Cambridge, 9 Sept. - 20 Sept. 2002.
The author is grateful to the lecturers of the school for exciting lectures and talks, as well as Isaac Newton Math. Institute for hospitality and financial support.
I am also grateful to the NATO ASI Conference on Bifurcations  of Periodic Orbits in Montreal,  Mathematical Institute of Stockholm University, Mittag-Lefler Institute (Druchholm), Uppsala University, CIRM (Lumini), KITPC/ITP-CAS,  Beijing,  Conference “Topological recursion and quantum algebraic geometry”, Aarhus Universitet for their support.
The author deeply grateful to  corresponding member of NASU S. Maksymenko for invitations to present talk on the Conference of the memory of V. Sharko and to present the corresponding text.
This paper has finished in Bulgaria. The author is deeply grateful to the Bulgarian Academy of Sciences, the Institute of Mathematics and Informatics of the Bulgarian Academy of Sciences, Professor P. Boyvalenkov for their support. 
The author was supported by Simons grant 992227.

\section{Braids, links and functions on manifolds}
\label{blfm}
Below in Sections 2, 3, 4 we follow to \cite{KasselTuraevl2008,KasselReuten2007,Artin1926,Markov1960,Markov1945,Turaev2002} and references therein. 
\begin{definition} (Configuration spaces of the ordered sets of points).
	Let $M$ be a closed smooth manifold	 and let $M^n$ be the product of $n$ spaces $M$ with the topology of the product. 
	Put 
	\begin{equation}
   \label{mpb}
	{\mathcal F}_n(M) = \{(u_1, \ldots , u_n) \in M^n; u_i \ne u_j, i \ne j \}.
	\end{equation}   
\end{definition}

\begin{remark}
	If $M$ is a topological space of the dimension $dim \; M$ (possibly with the boundary $\partial M$) then
	the dimension of ${\mathcal F}_n(M)$ is equal $n\cdot dim \; M$. The topological space ${\mathcal F}_n(M)$ is connected.
\end{remark}
\begin{definition}
  The fundamental group $\pi({\mathcal F}_n(M))$ of the manifold ${\mathcal F}_n(M)$ is called the group of pure braids with 
  $n$ strands.
\end{definition}

Let now $M$ be a connected topological manifolds of the dimension $\ge 2$, $M^{in} = M \setminus \partial M$, 
$Q_m \subset M^{in}$, $Q_m$ contains $m \ge 0$ points. Put 
\begin{equation*}
	{\mathcal F}_{m,n}(M) = {\mathcal F}_n(M \setminus Q_m).
	\end{equation*}   
	and for symmetric group $S_n$ put
\begin{equation*}
	{\mathcal G}_{m,n} = {\mathcal F}_{m,n}(M)/S_n .
	\end{equation*}   
\begin{definition}
  The fundamental group $\pi({\mathcal G}_{m,n})$  is called the  braids group of the manifold $M \setminus Q_m$ with 
  $n$ strands.
\end{definition}
\begin{remark}
	In the case of $M = {\mathbb R}^2$ we obtain groups of pure braids  and braids in the sense of E. Artin and A. Markov.
\end{remark}
Let $M$ be a three dimensional  topological manifold possible with the boundary $ \partial M$.
Recall that a geometric link in $M$ is a locally flat closed one dimensional submanifold in  $ M$.

The connections between algebraic functions and braids were also studied in the works
 of V. Arnol’d \cite{Arnold70} and others.

\section{Artin braids  (\cite{Artin1926,Markov1945,KasselTuraevl2008}) }

\begin{definition}
The group of braids of Artin $B_n$ is a group given by $n-1$ generators of $\sigma_1, \ldots \sigma_{n-1}$ satisfying the relations
\begin{equation*}
	\sigma_i \sigma_j =  \sigma_j \sigma_i 
	\end{equation*}   
	for integer $i, j = 1,\ldots, n-1, \; |i - j| > 1$ and 
	\begin{equation*}
	\sigma_i \sigma_{i+1} \sigma_i =  \sigma_{i+1} \sigma_i \sigma_{i+1}
	\end{equation*}   
	for $i = 1, \ldots , n-2$.
\end{definition}

\subsection{Generators and relations of groups}  Here and below we follow to  \cite{Markov1945,CoxeterMoser1972,MagKarSol1974,Buch1982,Koepf2021}.
Let $O = \{g_1, \ldots , g_m\}$ be a set of  generators of a discrete group $G$. Words are defined by concatenating symbols $g_i$ or  
$g_{i}^{-1}, \; i = 1, \ldots , m$. 
Let $ r( g_1, \ldots , g_m)$ be a word and $e$ be the unit of the group. \\
Let relations
\begin{equation*}
R = \{r_{1}( g_1, \ldots , g_m) = e, \ldots , r_{s}( g_1, \ldots , g_m) = e\}
	\end{equation*}   
and  generators  $O$ generate the group $G$. 
\begin{definition}
A pair $(O, R)$ of such $O$, $R$ is called the code of the group $G$.
\end{definition}
Recall well known
\begin{proposition}
	Let $ r_{1}( g_1, \ldots , g_m) , \ldots , r_{s}( g_1, \ldots , g_m)$ be words which represent relations and $F(O)$ be the free group which is generated by elements of $O$. Let $N$ be the smallest normal subgroup of $F(O)$ which contains words  $ r_{1}( g_1, \ldots , g_m) , \ldots , r_{s}( g_1, \ldots , g_m)$. Than the group $G$ is isomorphic to the factor group $F(O)/N$.
\end{proposition}
\subsection{The word problem for finitely presented groups}
Put
\begin{equation*} 
<O,R> = F(O)/N.
\end{equation*}  

\begin{problem}
	Given the code of the group $G$ and words
	 $u, \; v \in F(O)$.
	Decide whether u, v are equal in <O,R>.
\end{problem}

\subsection{Problems of  reduction and  simplification}

The problem of reduction and its solutions is based on the theory of reduction relations.
In most cases a reduction relation on a set is the binary relation on the set.
 \subsection{The word problem for braid groups}
 There are a number of solutions for the word problem in the braid groups by Garside, Thurston, and 
 others \cite{Garside1969,Epstein1998}.
\subsection{Braids as  cobordisms}
  Geometric interpretation of Artin braids give examples of cobordisms. Let $I = [0, 1]$ be the unit interval.
By the definition of cobordism 
, for a given element $b = \{b_i\}_{i = 1}^n$ (a braid on $n$ strands) of the braid group $B_n, n \ge 1$ we have the one dimensional compact manifold $W$ as the disjoint union of $n$ intervals $I$ and have two zero dimensional closed manifolds $M$ and $N$ consisting, respectively,  of $n$ points $b_i(0) = (i,0,1)$, $b_i(1) = (r(i),0,0)$, where $r(i)$ is some permutation. 
The tautological cobordism corresponding to a braid is constructed according to the structure of the braid and the definition of a cobordism,
 
 \section{Reidemeister, Alexander and Markov's results}

\subsection{Isotopy and isotopy classes}

\begin{definition} (Milnor)
	Let $f, \; g$ be two diffeomorphisms of manifolds: 
	\begin{equation*}
	f, g; X \to Y .
	\end{equation*}   
	A diffeomorphism $f$ is smoothly isotopic to  diffeomorphism $g$ if there exists a smooth homotopy 
	\begin{equation*}
	 F: X \times [0, 1] \to Y
	\end{equation*}   
	 between $f$ and $g$ such that for every $t \in [0, 1]$  there exists a map
	 \begin{equation*}
	x \to  F(t,x)
	\end{equation*}  
	 that maps $X$ homeomorphically onto $Y$.	
\end{definition}
Let $Ist(X,Y)$ be the set of  smoothly isotopic  diffeomorphisms from $X$ to $Y$.
\begin{remark}
	An isotopy is an equivalence relation, and thus the set $Ist(X,Y)$  is split into equivalence classes.
\end{remark}

\subsection{Braid diagrams}
Let $I$ be the closed interval $[0, 1]$.
\begin{definition}
	A braid diagram of $n$ strands is a set $\mathcal D \subset {\mathbb R}\times I$ split in the form of a union of $n$ topological
segments, called the threads of the diagram $\mathcal D$, for which the following conditions are satisfied:\\
1) the projection ${\mathbb R}\times I \to I$ maps each strand homeomorphically onto $I$;\\
2) each point of the product $\{1, \ldots n\} \times \{0, 1\}$ is the end point of some unique thread; \\
3) each point of the product ${\mathbb R}\times I$ belongs to at most two strands. At every point of intersection, these threads
intersect transversally.
\end{definition}
For braid diagrams Reidemeister moves  are defined, which were  defined by
  Reidemeister \cite{Reidemeister1932} for knots and knot diagrams.
\begin{theorem}(Reidemeister \cite{Reidemeister1932})
Two regular representations correspond to isotopic links if and only if the diagrams are related by isotopy (fixing the crossing points)  and by a finite sequence of the three Reidemeister moves.
\end{theorem}

 Alexander have presented the relation between links and closed braids.
\begin{theorem}( Alexander \cite{Alexander1923})
 Any link is  isotopic to a closed braid on some number of strands.
\end{theorem}

A.A. Markov  \cite{Markov1945} gave the description of the set of isotopy classes of oriented links in ${\mathbb R}^3$ in terms of braids. For manifolds of the dimension grater than 3 A.A. Markov  \cite{Markov1960} has proved the undecidability 
 of the problem of homeomorphy (see also   \cite{Markov1936,Morton1986,Traczyk1995,Vogel1990}).
 
 \section{On monodromy representation of braid groups and on some connections}
Here we present connections which correspond to  monodromy representation of braid groups as at the beginning of the Section (\ref{blfm}). For motivations and general ideas on Yang-Baxter and Knizhnik-Zamolodchikov connections and equations see \cite{Kohno1987},   \cite{Drinfeld1989,KasselReuten2007}   and references therein.

\subsection{On monodromy representation of braid groups}

Let $E \to X$ be the fibre bundle of a topological space $X$ with the fundamental group 
$\pi (X, x_0)$ at the point $x_0 \in X$ and with the fibre $F_{x_0}$ at the same point.

\begin{proposition}
If the connection $\nabla$ is flat (i.e. $K = 0$) then there exists the monodromy representation 
 $\pi (X, x_0) \to Aut(F_{x_0})$ which induced by the acting of the fundamental group 
$\pi (X, x_0)$  on the fibre $F_x$ at the point $x_0$.
\end{proposition}

\subsection{On infinitesimal braid group relations, and connections}
\label{ibgc}
For motivations and general ideas on Yang-Baxter and Knizhnik-Zamolodchikov equations see \cite{Kohno1987}    and \cite{Drinfeld1989}.
Below we are using the model of the configuration space (\ref{mpb}) of the form 
\begin{equation*}	
 {\mathcal F}_n = {\mathcal F}_n(M) = \{(z_1,\ldots z_n) \in {\mathbb C}^n| z_i \ne z_j,  i \ne j  \} 
 \end{equation*}
which regarded as a smooth manifold.  
 At first recall some known facts about connections on smooth algebraic varieties.We consider connections with de Rham differentials.
\begin{remark}
In Section (\ref{ccss}) we consider, follow to methods and results by B. Dwork, P. Griffths, A. Grothendieck, P. Deligne, Yu. Manin, N. Katz a more general context.
\end{remark}
\begin{definition}
For a one-dimensional sheaf of germs of differentials $\Omega^1_{S/k}$ and a coherent sheaf ${\mathcal F}$ on $S$, a connection on the sheaf ${\mathcal F}$ is a homomorphism of sheaves 
$$\Delta: {\mathcal F} \to \Omega^1_{S/k}\otimes {\mathcal F}   $$
such that if $f \in \Gamma(U,{\mathcal O}_S), g \in \Gamma(U,{\mathcal F})$, then $\Delta(fg) = f\Delta(g)+df\otimes g$.
\end{definition}

Let $\Omega^{i}_{S/k}$ be the sheaf of germs of $i-$differentials,
$$ \nabla^{i}(\alpha \otimes f) = d\alpha \otimes f + (-1)^{i}\alpha
\wedge \nabla(f). $$
   Then $ \nabla, \; \nabla^{i}$ define the sequence of homomorphisms:
\begin{equation}
 \label{CC}
 {\mathcal F} \rightarrow \Omega^{1}_{S/k} \otimes {\mathcal F}
 \rightarrow \Omega^{2}_{S/k} \otimes {\mathcal F} \rightarrow \cdots     ,.
\end{equation}
The map
$K = \nabla \circ \nabla^{1}: {\mathcal F}  \rightarrow \Omega^{2}_{S/k} \otimes {\mathcal F}$ is called
the curvature of the connection $\nabla.$
\begin{definition}
 A connection $\nabla$ is called the flat if the curvature of the connection is equal $0$.
\end{definition}

Below we consider the $n$-dimensional case.
For the complex vector space $W$ let $\{A_{i,j}\},1 \le  i < j \le n$ be automorphisms of $W$ which satisfy conditions
\begin{equation}
\label{LA}	
[A_{i,j}, A_{k,l}] = 0, \;
[A_{i,j}, A_{i,k} +  A_{j,k}] = [A_{j,k}, A_{i,j} +  A_{i,k}].
 \end{equation}
\begin{remark}
Below we assume that these conditions are satisfied.
 \end{remark}
Consider the connection
\begin{equation*}
 \nabla = d - \Gamma,
 \end{equation*}
on the trivial bundle 
\begin{equation*}	
 {\mathcal F}_n \times W,
 \end{equation*}
where
\begin{equation*}	
d\omega = \sum_{1 \le i < j \le n} \frac{A_{i,j}}{z_1 - z_j} (dz_i - dz_j)\omega,
 \end{equation*}
\begin{equation*}	
\Gamma = \sum_{1 \le i < j \le n} \frac{A_{i,j}}{z_1 - z_j} (dz_i - dz_j),
 \end{equation*}
\begin{proposition}
 The connection  $\nabla = d - \Gamma$  is flat. 
\end{proposition}

\begin{problem}
\label{pbgf}
Matrices with conditions \ref{LA} define for a given $n$ the Lie algebra $\mathfrak p$. By the Lie algebra it is possible to construct the local Lie group $\mathfrak P$. The expansion in group identity of the group law of the local Lie group $\mathfrak P$ gives the formal group law. The algebras, groups, and formal Lie groups mentioned above correspond to pure braid groups. What formal group laws correspond to pure braid groups?
\end{problem}

In the case of the braid group, we must consider the factormanifold 
\begin{equation*}
 {\mathcal G}_n = {\mathcal F}_n / S_n 
 \end{equation*}
where $S_n$ is the permutation group. Let the group $S_n$ acts from the left on the vector space $W$. Then there is the right action of $S_n$ on the trivial bundle ${\mathcal F}_n \times W$ by the formula
\begin{equation*}	 
(z_1,\ldots z_n), w)\sigma  = (z_{\sigma(1)},\ldots z_{\sigma(n)}), \sigma^{-1}(w)), \sigma \in S_n.
 \end{equation*}
If $d\omega$ is invariant under the action of $S_n$ then by the (Grothendieck) descend theory the connection 
$\nabla = d - \Gamma$ descends on the non-trivial vector bundle $({\mathcal F}_n \times W) / S_n$.

\begin{problem}
 This problem is similar to the Problem  (\ref{pbgf}).
	 Now algebras, groups, and formal Lie groups  correspond to  braid groups. What formal group laws correspond to  braid groups in the monodromy representations?
	\end{problem}

\section{Connections on coherent sheaves of smooth schemes}
\label{ccss}
Here we follow to methods and results of B. Dwork, A. Grothendieck, P. Deligne, Yu. Manin, N. Katz \cite{Katz2007}, as well as the works \cite{KedlayaTuitman2012,Glazunov2002}.
In \cite{Glazunov2002} we considered applications of the methods of algebraic geometry to string theory.
An open string can be   interpreted as a unit (or trivial) braid group $B_1$.
A closed string is the closure of the element of the unit braid group.
Introduction to the consideration of non-trivial braid groups extends and deepens
emerging algebraic picture.
Below we use notations and results from subsection (\ref{ibgc}).

Let  $\Theta^1_{S/k}$  be the tangent sheaf on the scheme $S$, $\partial \in  \Gamma(U,\Theta^1_{S/k})$.
In this case, we can give a dual definition of a connection on a sheaf.
\begin{definition}
In the above notation, a connection is a homomorphism 
$$ \rho: \Theta^1_{S/k} \to End_{{\mathcal O}_S} ({\mathcal F}, {\mathcal F})  $$
such that $\rho(\partial)(f g) = \partial (f) g + f \rho(\partial)$.
\end{definition}
Recall that the cochain complex
$$ (K^{\bullet}, d) = \{ K^{0} \xrightarrow{d} K^{1} d = 0\xrightarrow{d} K^{2} \xrightarrow{d} \cdots \} $$
in the category of Abelian groups is called the sequence of Abelian groups and morphisms
$d: K^{p} \xrightarrow{d} K^{p+1}$ such that $d\circ d =0$.
\begin{remark}

Let ${\mathcal A}$ be an abelian
category, ${\mathcal K}({\mathcal A})$ the category of complexes over
${\mathcal A}$. Furthermore, there are various full subcategories of
${\mathcal K}({\mathcal A})$ whose respective objects are the complexes
which are bounded below, bounded above, bounded in both sides. The
notions of homotopy morphism, homotopical category, triangulated
category are described in \cite{GeM88}. The bounded derived
category ${\mathcal D}^{b}(X)$ of coherent sheaves on $X$ has the
structure of a triangulated category \cite{GeM88}.
\end{remark}
\subsection{Integration of Connections}
\begin{definition}
  A connection is {\em integrable} if (\ref{CC}) is a complex.
\end{definition}
\begin{proposition}
\label{IC}
 The statements a),b), c) are equivalent: \\
a)  the connection $\nabla$ is  integrable;\\
b)  $K = \nabla \circ \nabla^{1} = 0;$\\
c)  $ \rho $ is the Lie-algebra homomorphism of sheaves of Lie algebras.\\
\end{proposition}
 More generally, let $X, S$ be smooth schemas over $k, f: X \rightarrow S $
a smooth morphism,
\begin{equation}
 \label{C1}
  \Omega_{X/S}^{\bullet}: {\mathcal O}_{X} \rightarrow \Omega^{1}_{X /S}
 \rightarrow \Omega^{2}_{X /S}  \rightarrow \cdots     ,
\end{equation}
the de Rham complex of relative differentials.
There is the exact sequence
\begin{equation}
 \label{C2}
  0 \rightarrow f^{*}(\Omega_{S}^{\bullet})  \rightarrow \Omega_{X}^{\bullet}
 \rightarrow \Omega_{X/S}^{\bullet} \rightarrow 0.
\end{equation}
In the case there is an integrable connection, the Gauss-Manin connection. \\

By an abelian sheaf we mean a sheaf of abelian groups.
Let $R^{0}f_{*}$ be the functor from the category of complexes
of abelian sheaves on $X$ to the category of abelian sheaves on $S.$ Denote by
 ${\mathcal H}^{i}_{DR}(X/S)$ the sheaf of de Rham cohomologies such that
$${\mathcal H}^{i}_{DR}(X/S) = R^{i} f_{*}\Omega_{X/S}^{\bullet} .$$
Recall that
$${\mathcal H}^{i} = {\mathcal H}^{i}_{DR}= R^{i} f_{*}\Omega_{X/S}^{\bullet} $$
is called the Gauss-Manin bundle.

Here $ R^{i} f_{*}$   are the hyperderved functor of  $R^{0}f_{*}.$

For each $i\ge 0,{\mathcal H}^{i}_{DR}( X/S )$ is a locally free coherent algebraic sheaf on
$S,$ whose fiber at each point $s \in S$ is the ${\bf C}-$vector  space
$H^{i}(X_{s}, {\bf C})$ and has the Gauss-Manin connection.
${\mathcal H}^{i}_{DR}( X/S )$ has the main interpretation as the Picard-Fuchs equations and
$H^{i}(X_{s}, {\bf C})$ in the  interpretation is the local system of germs of solutions of the equations.

There is the canonical filtration of  $\Omega_{X/S}^{\bullet} $ by locally free subsheaves
$$\Omega_{X/S}^{\bullet}  = F^{0}(\Omega_{X/S}^{\bullet}) \supset   F^{1}(\Omega_{X/S}^{\bullet}) \supset \cdots $$
given by
$$F^{i}(\Omega_{X}^{n}) := Im(f^{*}(\Omega_{S}^{i}) \otimes \Omega_{X}^{n-i} \rightarrow \Omega_{X}^{n}). $$

\section{Formal groups in  complex cobordism}

\subsection{Algorithm for constructing of Lazard`s one dimensional universal formal group}
Here we follow to \cite{Lazard1955,Hazewinkel1978,Glazunov2018}.
We extract from Lazard \cite{Lazard1955} an algorithm of constructing of Lazard`s one dimensional universal formal group law.
But we give here the structute of the Lazard algorithm only.
 The constructions of $n (n\ge 1)$-buds, formal groupoids and moduli spaces of one dimensional formal groups are  investigated and are used.
 Let $A_q$ and $A^{`}_q$ be the rings of polynomials $A_q = {\mathbb Z}[\alpha_1, \ldots, \alpha_q]$ and $A^{`}_q = {\mathbb Q}[\alpha_1, \ldots, \alpha_q]$.
 \begin{proposition}
	The structure of the algorithm has the next form. For the (one-dimensional) $1$-bud   $x + y + \alpha_1 xy$ we put: 
\begin{equation*}	
 f_1(x,y) = x + y + \alpha_1 xy; \;  
  \varphi_1 = x - \frac{1}{2}\alpha_1 x^2; 
 \end{equation*}
 \begin{equation*}
  f_{q + 1}(x,y) = f_q(x,y) + h^{`}(x,y) + \alpha_{q+1}C_{q + 2}(x,y), \; f_q(x,y) \in A_q;
 \end{equation*}
 compute $ \varphi_{q + 1},  \;  \varphi_{q + 1} \in A^{`}_q$.
\end{proposition}
\begin{definition}
 The ring ${\mathbb L} = {\mathbb Z}[\alpha_1, \alpha_2,  \cdots ]$ is called the  Lazard ring.
\end{definition}
\subsection{On formal groups in complex cobordism}
Let $M$ be a smooth manifoldd, $TM$ be the tangent bundle on $M$ and ${\mathbb R}^m$ the trivial real $m$-dimensional bundle of $M$.
\begin{definition}
	A manifold $M$ is \emph{stably complex} if  for some natural $m$ the real vector bundle $TM\bigoplus{\mathbb R}^m$ admits a complex structure.
\end{definition}
Let $M_1$ and $M_2$ be two smooth manifolds of dimension $n$, and $W$ be the smooth manifold of dimension $n+1$ with a boundary that is the union of $M_1$ and $M_2$, i.e. $\partial W  = M_1 + M_2$.
\begin{definition}
	Let  $M_1, M_2, W$ be stably complex manifolds. In above notations the \emph{complex (unitary) cobordism} between $M_1$ and $M_2$ is a manifold $W$ whose boundary is the disjoint union of $M_1$, $M_2$, $\partial W  = M_1 +
	{\overline M_2}$ where the corresponding structure on  $\partial W$ is induced from $W$ and ${\overline M_2}$ denotes the manifold with opposite structure.
\end{definition}
\begin{remark}
	Suppose we have the relation of complex cobordism. Then the relation divides stably complex manifolds on equivalence classes called classes of complex cobordisms.
\end{remark}
\begin{lemma}
	The set of classes of complex cobordisms with operations of disjoint union and product of stable manifolds form commutative graded ring  ${\mathbb L} = {\mathbb Z}[v_1, v_2,  \cdots ]$ .
\end{lemma}
\begin{theorem}
(A.S.  Mishchenko) Let $g(t)$ be the logarithm of Lazard universal  formal group, and $[{\mathbb C}P^n] $ are classes of unitary  cobordisms of complex projective spaces.
Then
\begin{equation*}	
 g(t) = \sum_{n \ge 0} \frac{[{\mathbb C}P^n]}{n+1} t^{n+1}, \;
 [{\mathbb C}P^0] = 1.
\end{equation*}
\end{theorem}

\section*{Formal groups in  algebraic cobordism}
Here we follow to \cite{LevineMorel2007}.
Let $k$ be a field of characteristic zero and ${\mathbf {Sm}}(k)$ be the full subcategory of smooth quasi-projective $k$-schemes of the category of separable finite-type $k$-schemes.
Let $A^*$ be an oriented cohomology theory on ${\mathbf {Sm}}(k)$  and let $c_1(L)$ be the first Chern class of line bundle 
$L$ on    $X \in {\mathbf {Sm}}(k)$.

 \begin{proposition}
( Quillen, \cite{LevineMorel2007}).
Let $L,\; M$ be   line bundles on $X \in {\mathbf {Sm}}(k)$.. There exists formal group law $F_A$ with coefficients in 
 $A^*$ such that
\begin{equation*}	
c_1(L\otimes M) = F_A(c_1(L),c_1(M)).
 \end{equation*}
\end{proposition}

\begin{theorem}
(Levine-Morel).There is a universal oriented cohomology theory $\Omega$ over $k$ called algebraic cobordism.
The classifying map $\phi_{\Omega}: {\mathbb L} \to \Omega^*(k)$ is an isomorphism , so $F_{\Omega}$ is the universal formal group law.
\end{theorem}

\begin{problem}
	It seems that very little is known about applying  $n (n \ge 2)$-dimensional commutative formal groups to 
	cobordism theory. For instance what is the application of a two-dimensional 1-bud ?
	\end{problem}




\bibliography{mybib}{}
\bibliographystyle{plain}


\begin{center}
{\small {\rm Glushkov Institute of Cybernetics National Academy of Sciences of Ukraine,\\
Institute of Mathematics and Informatics Bulgarian Academy of Sciences }}
\end{center} 

\end{document}